\documentclass{article}
\usepackage{spconf,amsmath,graphicx,amsfonts}

\title{Precise Error Analysis of the $\ell_2$-LASSO}
\vspace{-10pt}
\name{ Christos Thrampoulidis$^*$
Ashkan Panahi$^\dagger$, Daniel Guo$^*$, Babak Hassibi$^*$ \thanks{The work of B. Hassibi was supported in part by the National Science Foundation under grants CNS-0932428, CCF-1018927, CCF-1423663 and CCF-1409204, by the Office of Naval Research under the MURI grant N00014-08--0747, by the Jet Propulsion Lab under grant IA100076, by a grant from Qualcomm Inc., and by King Abdulaziz University.}
}
\vspace{-10pt}
\address{
* Department of Electrical Engeeniring, Caltech, Pasadena, USA\\
$\dagger$ Signal Processing Group, Chalmers University of Technology, Gothenburg, Sweden
}
\usepackage{enumerate}
\usepackage{amsmath,amssymb,algorithm,cite,cases,bm,float,graphicx,url}
\usepackage{pbox}

\usepackage{amsthm}
\newtheorem*{theorem*}{Theorem}

\newtheorem*{lemma*}{Lemma}

\newtheorem{defn}{Definition}[section]

\newtheorem{thm}{Theorem}[section]

\usepackage{lipsum}
\usepackage{amsmath,amssymb,algorithm,cite,caption,cases,bm,float,graphicx,url,color}
\usepackage[hidelinks]{hyperref}
\usepackage{comment}
\usepackage{thm-restate, enumerate }

\definecolor{darkred}{RGB}{150,0,0}
\definecolor{darkgreen}{RGB}{0,150,0}
\definecolor{darkblue}{RGB}{0,0,200}
\hypersetup{colorlinks=true, linkcolor=darkred, citecolor=darkgreen, urlcolor=darkblue}

\usepackage{graphicx}
\usepackage{caption}
\usepackage{subcaption}

\usepackage[numbers]{natbib}

\newcommand{\beq}{\begin{equation}}
\newcommand{\eeq}{\end{equation}}
\newcommand{\bea}{\begin{align}}
\newcommand{\eea}{{\end{align}}}

\newcommand{\lah}{\hat{\la}}

\newcommand{\Db}{{\mathbf{D}}}

\newcommand{\alphas}{\alpha_*}

\newcommand{\betas}{\beta_*}

\newcommand{\Ec}{\mathcal{E}}
\newcommand{\Rc}{\mathcal{R}}
\newcommand{\Gco}{\mathcal{\Gc}_o}

\newcommand{\ws}{{\w}_*}

\newcommand{\vp}{\vspace{-5pt}}

\newcommand{\Dlf}{{\mathbf{D}}(\la )}

\newcommand{\Clf}{{\mathbf{C}}(\la )}

\newcommand{\lac}{\lambda_{\text{crit}}}
\newcommand{\labb}{\lambda_{\text{best}}}

\newcommand{\A}{\mathbf{A}}

\newcommand{\w}{\mathbf{w}}

\newcommand{\x}{\mathbf{x}}
\newcommand{\ub}{\mathbf{u}}
\newcommand{\g}{\mathbf{g}}
\newcommand{\vb}{\mathbf{v}}

\newcommand{\y}{\mathbf{y}}

\newcommand{\z}{\mathbf{z}}

\newcommand{\h}{\mathbf{h}}

\newcommand{\Dc}{\mathcal{D}}
\newcommand{\Lc}{\Phi}

\newcommand{\Nn}{\mathcal{N}}

\newcommand{\Gc}{\phi}

\newcommand{\R}{\mathbb{R}}
\newcommand{\Pro}{{\mathbb{P}}}
\newcommand{\E}{{\mathbb{E}}}

\newcommand{\la}{{\lambda}}

\newcommand{\eps}{\epsilon}

\newcommand{\Fc}{\mathcal{F}}

\newcommand{\nn}{\nonumber}

\begin{document}

%

\maketitle
\begin{abstract} 
A classical problem that arises in numerous signal processing applications
asks for the reconstruction of an unknown, $k$-sparse signal $\x_0\in \R^n$ from underdetermined, noisy, linear measurements $\y=\A\x_0+\z\in \R^m$. One standard approach is to solve the following convex program $\hat\x=\arg\min_\x \|\y-\A\x\|_2 + \la \|\x\|_1$, which is  known as the $\ell_2$-LASSO. We assume that the entries of the  sensing matrix $\A$ and of the noise vector $\z$ are i.i.d Gaussian with variances $1/m$ and $\sigma^2$. In the large system limit when the problem dimensions grow to infinity, but in constant rates, we \emph{precisely} characterize the limiting behavior of the normalized squared error $\|\hat\x-\x_0\|^2_2/\sigma^2$. Our numerical illustrations validate our theoretical predictions.


\end{abstract}
 

\vspace{-10pt}
\section{Introduction}

\vspace{-15pt}
\subsection{Motivation}
The Least Absolute Shrinkage and Selection Operator (LASSO) is a celebrated convex progam used to estimate sparse signals from noisy linear underdetermined observations. Given a vector of observations $\y = \A\x_0 + \z \in \R^m$ of an unknown, but $k$-\emph{sparse} (i.e., at most $k$ nonzero entries), signal 
$\x_0\in\R^{n}$,
 the $\ell_2$-LASSO\footnote{Also known as ``square-root LASSO" \cite{Belloni}. Please refer to \cite[Sec.~1.3]{OTH}.} produces the following estimate for $\x_0$:
\begin{align}\label{eq:LASSO}
\hspace{-10pt}\hat\x
 := \arg\min_{\x} \Lc(\x;\A,\z) :=\|\y - \A\x\|_2 + \frac{\la}{\sqrt{m}} \|\x\|_1.
\end{align}
Here, $\A\in\R^{m\times n}$ is the sensing matrix, $\z\in\R^m$ is the noise vector and $\la> 0$ is a regularizer parameter.  The LASSO has been long investigated from different perspectives and shown to enjoy unique properties, in terms of both computation and precision. Yet, some important asymptotic properties of it have not yet been fully understood. Our interest is on the \emph{exact} characterization of 
the reconstruction error $\|\hat\x-\x_0\|_2$. 


\vp
\subsection{Contribution}
\vp
We assume a generic setup in which the entries of the sensing matrix and the \emph{non-zero} entries of $\x_0$  are i.i.d  Gaussian. Also, the noise vector $\z$ is assumed to have entries i.i.d Gaussian with variance $\sigma^2$.
Under this assumption we derive an \emph{asymptotically exact} expression for the normalized squared error (NSE) $\|\hat\x-\x_0\|_2^2/\sigma^2$ of the $\ell_2$L-ASSO. 
In the low noise regime $\sigma^2\rightarrow 0$, our result reduces to simple and interpretable formulae.
Although our theoretical analysis requires an asymptotic setting in which the problem dimensions grow to infinity, our numerical illustrations suggest that the predictions be accurate already with problem dimensions ranging over only a few hundreds.  Also, we remark on our assumption on the Gaussian nature of the sensing matrix. This assumption has a long tradition in the statistics literature and sheds important insights \cite{MontanariLASSO}; the Gaussian ensemble has wonderful properties which make the analysis tractable, while at the same time many of the results 
generalize to a wider class of distributions.
The main technical tool used in our analysis is a gaussian comparison inequality due to Gordon \cite{Gor}. When combined with appropriate convexity assumptions, the inequality can be shown to be tight \cite{stojnic2013meshes,OTH,Allerton,tight}, which makes it ideal for our precise analysis.



\vspace{-10pt}
\subsection{Relevant  Literature}
\vp
The LASSO was introduced by Tibshirani in \cite{TibLASSO} in the form
\vp
\begin{align}
\hat\x = \arg\min_{\x} \|\y-\A\x\|_2 \text{ s.t. } \|\x\|_1\leq\|\x_0\|_1.\label{eq:conLASSO}
\end{align}
\vp
\eqref{eq:LASSO} is a regularized version of \eqref{eq:conLASSO}. An alternative to \eqref{eq:LASSO} solves  
\vspace{-8pt}
\begin{align}
\hat\x = \arg\min_{\x} \|\y-\A\x\|_2^2+ \tau\|\x\|_1.\label{eq:ell22LASSO}
\end{align}
Both this and \eqref{eq:LASSO} are variations of the same algorithm and Lagrange duality ensures that they both become equivalent
to the constrained optimization \eqref{eq:conLASSO} for proper choice of the regularizer parameters.
There are reasons to argue in favor of either of them \cite{Belloni,TOH14}, but they go beyond the scope of this paper. 
Early well-known bounds on the reconstruction error of the LASSO were  order-wise in nature (i.e. accurate only up to constant multiplicative factors) and derived based on RIP and Restricted Eigenvalue assumptions on the measurement matrix \cite{candes2006stable,candes2007dantzig,bickel,negahban2012unified,Belloni}. To the best of our knowledge, the first precise asymptotic evaluation of the limiting behavior of the LASSO reconstruction error is due to Bayati and Montanari \cite{MontanariLASSO}; they consider problem \eqref{eq:ell22LASSO}, i.i.d Gaussian sensing matrix $\A$ and use the Approximate Message Passing (AMP) framework. More recently, and closer to our work, Stojnic introduced an alternative framework that cleverly uses a celebrated gaussian comparison inequality due to Gordon \cite{Gor}; he applied this on \eqref{eq:conLASSO} and obtained precise characterizations for its worst-case NSE, which he showed to occur in the limit $\sigma^2\rightarrow 0$. The same framework was used in \cite{OTH} to generalize the results of \cite{StoLASSO} to arbitrary convex regularizer functions. Moreover, the authors in \cite{OTH} consider \eqref{eq:LASSO} and obtain simple, precise error formulae when $\sigma^2\rightarrow 0$. Our work, extends this result of \cite{OTH} in several directions for the case of sparse recovery. First, it holds for arbitrary values of the noise-variance $\sigma^2$; to our knowledge, this is the first such precise asymptotic result for \eqref{eq:LASSO} (\cite{MontanariLASSO} considers \eqref{eq:ell22LASSO}), and suggests the capabilities of the framework used.  Also, when $\sigma^2\rightarrow 0$, it recovers the result of \cite{OTH} and extends it on the range of values of $\la$ for which it holds. Finally, we believe that our result can be used to prove that the worst-case NSE of \eqref{eq:LASSO} occurs when $\sigma^2\rightarrow 0$. This would imply that the simple and interpretable expressions corresponding to that regime are tight bounds on the NSE for arbitrary noise levels. Some additional effort is required to prove this claim and is, thus, left for future work.

\vspace{-18pt}
\section{Analysis}

\vspace{-10pt}
\subsection{Problem Setup}\label{sec:setup}
\vspace{-5pt}
For the rest of the paper, let $\Nn(\mu,\sigma^2)$ denote the normal distribution of mean $\mu$ and variance $\sigma^2$. Also, we use $\|\cdot\|$ instead of $\|\cdot\|_2$. Let $\x_0\in\R^n$ denote the unknown signal and $\y=\A\x_0 + \z\in\R^m$ denote the vector of observations. We make the following assumptions:
\vspace{-6pt}
\begin{itemize}
\item the entries of $\A$ are i.i.d $\Nn(0,1/m)$,
\vp
\item the entries of $\z$ are i.i.d $\Nn(0,\sigma^2)$,
\vp
\item $\x_0$ is $k$-sparse with \emph{fixed} support set $S$, i.e. $(\x_0)_i=0$ for all $i\notin S$. Also, $(\x_0)_i\overset{\text{i.i.d}}{\sim}\Nn(0,1)$,
\end{itemize}
\vspace{-8pt}
Consider estimating $\x_0$ via the solution $\hat\x:=\hat\x(\la)$ of the LASSO program in \eqref{eq:LASSO}.
Our goal is to provide tight expressions for the reconstruction error $\|\hat\x-\x_0\|$;
this depends explicitly on $\x_0$ and implicitly on $\A$ and $\z$.
We consider an asymptotic setting in which the problem parameters $n,m$ and $k$ grow proportionally as $m/n\rightarrow\delta \in(0,1)$ \footnote{Our analysis extends to the overdetermined case, where $\delta\in[1,\infty)$. For simplicity, in this paper we focus on the  underdetermined regime. } and $k/m\rightarrow \gamma\in(0,1)$. Our result characterizes the limiting behavior of the quantity of interest as $n\rightarrow\infty$. To suppress notation, we use the symbol ``$\approx$" to denote convergence in probability, i.e. $X_n\approx X$ is used to denote that a sequence of random variables $X_n$ converges in $X$ in that for all $\eps>0$, $\lim_{n\rightarrow\infty}\Pro\left( |X_n - X|>\eps X \right) = 0$. Similarly, we write $X_n\gtrsim X$ if for all $\eps>0$, $\lim_{n\rightarrow\infty}\Pro\left( X_n < (1-\eps) X \right) = 0$.

\vspace{-12pt}
\subsection{Preliminaries}

\vp
It is convenient to rewrite \eqref{eq:LASSO} changing the optimization variable to the error vector $\w:=\x-\x_0$:
\footnote{Also, note the re-scaling in \eqref{eq:main} with a factor of $\sqrt{m}$. 
}:
\begin{align}\label{eq:main}
\min_{\w} \Lc(\w;\A,\z) := \sqrt{m}\|\z-\A\w\| + \la \|\x_0 + \w\|_1.
\end{align}
Denote $\hat\w^{\la,\sigma}(\A,\z)$ the solution of \eqref{eq:main}. We often drop the dependence on the arguments $\la,\sigma,\A,\z$ when clear from context.
Suppose we want to show that $\|\hat\w\|\approx \alphas$, for  appropriate $\alpha_*\geq 0$. Equivalently, for all $\eps>0$ and sets $\Rc_\eps:=\{\ell~|~ |\ell -\alphas| > \eps\alphas\}$, we wish that $\lim_{n\rightarrow\infty}\Pro( \|\hat\w\|\in\Rc_\eps ) = 0$. It suffices to prove the existence of $\delta:=\delta(\eps)>0$ such that  
\vspace{-10pt}
\begin{align}\label{eq:prob_state}
\lim_{n\rightarrow \infty}\Pro\left( \min_{\|\w\|\in \Rc_\eps} \Lc(\w) \leq (1+\delta ) \Lc(\hat\w) \right) = 0.
\end{align}
Directly showing \eqref{eq:prob_state} through analyzing $\Lc(\w)$ is difficult. Instead, we use the Gaussian min-max Theorem to translate the LASSO objective function in \eqref{eq:main} to a simpler one, that is amenable to direct analysis.

\vspace{-10pt}
\subsection{Introducing a Simpler Optimization}
\vp
The Gaussian min-max Theorem belongs to the family of the so called Gaussian comparison inequalities \cite[Ch.~3]{Ledoux} and was proved by Gordon in \cite[Lemma 3.1]{Gor}. To see how that theorem can be applied in our case, observe that
\vspace{-8pt}
$$
\Lc(\w) = \max_{\|\ub\|\le 1} \ub^T[\sqrt{m}\A, -\z/{\sigma}]\left [
\begin{array}{c}
\w \\
\sqrt{m}\sigma \\
\end{array}
\right ]+\lambda\|\x_0+\w\|_1.
$$
Let $\g\in\R^m$, $\h\in\R^n$, $g\in\R$ have i.i.d $\Nn(0,1)$ entries and
\vp
\begin{align}
&\Gc(\w;\g,\h) := \max_{\|\ub\|\le 1}\Big\{ \sqrt{\|\w\|^2+m\sigma^2}\g^T\ub \nn\\
&\qquad\qquad-\|\y\|(\h^T\w+\sqrt{m}\sigma g)+\lambda\|\w+\x_0\|_1 \Big\}.\label{eq:Gw}
\end{align}
Then, the Gaussian min-max Theorem\footnote{In fact, the statement in \eqref{eq:low} is a slight variation of the original statement of the Gaussian min-max Theorem. Please refer to \cite{tight} for details.} states that for arbitrary set $\Rc\subset\R^n$ and $c\in\R$:
\begin{align}\label{eq:low}
\Pro( \min_{\w\in\Rc} \Lc(\w;\A,\z) \leq c)\leq 2\Pro( \min_{\w\in\Rc} \Gc(\w;\g,\h) \leq c ).
\end{align}
 Recently, it was shown in \cite{StoLASSO,tight} that when combined with appropriate convexity assumptions, the Gaussian min-max Theorem is \emph{tight}. In particular, for a \emph{convex} compact set $\Rc$ it is also true \cite[Theorem II.1]{tight} that 
\vp
\begin{align}\label{eq:up}
\Pro( \min_{\w\in\Rc} \Lc(\w;\A,\z) \geq c)\leq 2\Pro( \min_{\w\in\Rc} \Gc(\w;\g,\h) \geq c ).
\end{align}
\eqref{eq:low} and \eqref{eq:up} are critical for establishing \eqref{eq:prob_state}. They suggest the analysis of the following optimization problem\footnote{ The proof of \eqref{eq:up} in \cite{tight} requires the set $\Rc$ to be compact. This technical detail can be resolved by assuming a sufficiently large upper bound $K$ on $\|\hat\w\|$ such that constraining the minimization in \eqref{eq:2ana} over the set $\|\w\|\leq K$ does not change the optimal cost.}, which we refer to as ``Gordon's optimization (GO)":
\vp
\begin{align}\label{eq:2ana}
(GO)\quad\quad\Gc(\ws;\g,\h):= \min_{\w}\Gc(\w;\g,\h),
\end{align}
in place of the original LASSO optimization in \eqref{eq:main}. To see how this can be useful, assume that $\Gc(\ws)\approx \Dc_*$; then, it would follow directly from \eqref{eq:low} and \eqref{eq:up} that $\Lc(\hat\w)\approx \Dc_*$. Next, we analyze such asymptotic properties of (GO) .

\begin{figure*}[!t]
\footnotesize
%
\begin{align}\label{eq:defn}
&\Dlf := k(1+\la^2) + (n-k)\rho(1,\la),\quad
\Clf := -(\la/2) \partial\Dlf/\partial\la, \quad\psi(\la) := \la/Q^{-1}\left( {1}/{2} + ({m-\Dlf-\Clf})/({2k}) \right),\nn\\
&f(\la,\phi) := m- \Dlf + m{\sigma^2}(1-\phi^2) + k\left(\phi^2 - \la^2- 2\right)\left(2Q\left({\la}/{\phi}\right)-1\right) +  k \sqrt{{2}/{\pi}}{\lambda\phi}\exp\left(-{\la^2}/({2\phi^2})\right). &
\end{align}
\hrulefill
\vspace*{4pt}
\end{figure*}

\vspace{-8pt}
\subsection{Analyzing Gordon's Optimization}
\vspace{-5pt}
\subsubsection{Scalarization}
\vspace{-5pt}
We begin with simplifying (GO) and reducing it to an optimization problem involving only scalars. For now, assume that $\g,\h$ and $g$ are all fixed. 
Note that the maximization over the direction of $\ub$ in \eqref{eq:Gw} is easy to evaluate. Furthermore, we may write $\|\x_0+\w\|_1 = \max_{\|\vb\|_\infty\leq 1}(\x_0+\w)^T\vb$. With these:\vp
\begin{align}\label{eq:Gw2}
\Gc(\w_*) &= \min_{\w}\max_{\substack{0\leq\beta\leq 1\\ \|\vb\|_\infty\leq 1}}\Big\{ \sqrt{\|\w\|^2+m\sigma^2}\|\g\|\beta \\ \vp
&\quad\quad-\beta(\h-\la\vb)^T\w-\beta\sqrt{m}\sigma g+\lambda\x_0^T\vb\Big\}.\nn
\end{align}
\vspace{-2pt}
The objective function in \eqref{eq:Gw2} is convex-concave in $\w$ and $\beta,\vb$. Also, the constraint sets are closed convex and one of them is bounded; thus, we can switch the correspoding order of min-max \cite[Cor.~
37.3.2]{Roc70}. After this, the minimization over the direction of $\w$ is easy to evaluate:
\vspace{-8pt}
\begin{align}
\Gc(\ws) &= \max_{\substack{0\leq\beta\leq 1\\ \|\vb\|_\infty\leq 1}}\min_{\alpha\geq 0} \Big\{ \sqrt{\alpha^2+m\sigma^2}\|\g\|_2\beta\nn-\beta\sqrt{m}g\sigma
\\ 
&\qquad\qquad\qquad-\alpha\|\lambda \vb-\beta\h\|_2+\lambda \x_0^T\vb\Big\}.\nn
\end{align}
Note that the optimization variable $\alpha$ above plays the role of the $\ell_2$-norm of $\w$; its optimal value $\alphas$ is equal to $\|\w_*\|_2$. Switching one last time the order of optimization between $\vb,\beta$ and $\alpha$, we conclude with the following optimization:
\begin{align*}
\min_{\alpha\geq 0}\max_{\substack{0\leq\beta\leq 1}}\left\{ \sqrt{\alpha^2+m\sigma^2}\|\g\|_2\beta-\beta\sqrt{m}g\sigma-\beta\Fc(\alpha,\beta)\right\},
\end{align*}
where $\Fc(\alpha,\beta) := \min_{\|\vb\|_\infty\le 1}\left\{\alpha\|\frac{\lambda}{\beta} \vb-\h\|_2-\frac{\lambda}{\beta} \x_0^T\vb\right\}.$ The final step amounts to reducing the optimization over $\vb$ to a scalar optimization. The idea here is to apply the fact that for any real $r\ge 0$, one has $\sqrt{r}=\min_{p>0}\left(\frac{p}{2}+\frac{r}{2p}\right)$ in order to make the objective function separable. With this:
\vp
\begin{align}
&\mathcal{F}(\alpha,\beta)=\min_{\|\vb\|_\infty\le 1}\min_{p>0}\frac{\alpha p}{2}+\frac{a\|\lambda/\beta \vb-\h\|_2^2}{2p}-\frac{\lambda}{\beta} \x_0^T\vb\nn\\
&=\min_{p>0}\bigg \{\frac{\alpha p}{2}+\sum_{i\in S}\Big[\frac{\alpha}{2p}\min_{|\vb_i|\le \frac{\lambda}{\beta}}\left(\vb_i-\h_i-\frac{p(\x_0)_i}{\alpha}\right)^2\nn
\\
&-\frac{p(\x_0)_i^2}{2\alpha}-\h_i(\x_0)_i\Big]+\frac{\alpha}{2p}\sum_{i\notin S}\min_{|\vb_i|\le \frac{\lambda}{\beta}}\left(\vb_i-\h_i\right)^2 \bigg \}.\label{eq:above}
\end{align}
The second equality above follows from standard completion of squares. The scalar minimizations over $\vb_i$'s in \eqref{eq:above} are simple soft-thresholding operations: $\min_{|q|\leq\tau}(r-q)^2 = 
( (|r|-\tau)^+ )^2.$
Combining all the above, we have shown that
$
\Gc(\ws) = \min_{\alpha\geq 0}\max_{\substack{0\leq\beta\leq 1, p>0}} \Gco(\alpha,\beta,p),
$
for an appropriately defined objective function $\Gco(\alpha,\beta,p)$.

\vspace{-10pt}
\subsubsection{Concentration}
\vp
Our next step is to analyze the limiting behavior of $\Gco(\alpha,\beta,p)$. Recall that $\g,\h,g$  have i.i.d $\Nn(0,1)$ entries. Thus, $\|\g\|\approx\sqrt{m}$ and $g\approx 0$. Using these and applying the Law of Large Numbers to the summations in \eqref{eq:above} it can be shown that
$\Gco(\alpha,\beta,p)\approx \Dc(\alpha,\beta,p),
$ where 
\begin{align}\nn
&\Dc(\alpha,\beta,p):=\beta\Big( \sqrt{\alpha^2+m\sigma^2}\sqrt{m} - \frac{\alpha p}{2} +  k\frac{p}{2\alpha} \nn\\
& - k\frac{\alpha}{2p} \rho(\sqrt{1+p^2/\alpha^2} , \la/\beta  ) -(n-k) \frac{\alpha}{2p} \rho( 1 , \la/\beta  ) \Big),\nn
\end{align}
\vspace{-10pt}
\begin{align}\label{eq:D}
\rho(c,\tau)& := \E_{g\sim\Nn(0,1)}\left[ (|cg|-\tau)^+\right]^2\nn\\
 &= 2(c^2+\tau^2)Q(\tau/c) - \sqrt{2/\pi} c\tau e^{-{\tau^2}/{(2c^2)}},
\end{align}
and $Q(x) = (1/\sqrt{2\pi})\int_{x}^{\infty}e^{-x^2/2}\mathrm{d}x$.
Albeit some technicalities involved (skipped due to space limitations):
\vspace{-8pt}
\begin{subequations}
\begin{align}
\min_{\w} \Gc(\w)&\approx\min_{\alpha\geq 0}\max_{\substack{0\leq\beta\leq 1\\ p>0}}\Dc(\alpha,\beta,p),\label{eq:f1}\\
\min_{\|\w\| \in \Rc_\eps} \Gc(\w)&\approx\min_{\alpha\in\Rc_\eps}\max_{\substack{0\leq\beta\leq 1\\ p>0}}\Dc(\alpha,\beta,p)\label{eq:f2}.
\end{align}
\end{subequations}
\vp
with $\Rc_\eps = \{\ell ~|~ |\ell -\alpha_*|>\eps\alpha_* \}$ and $\alphas$ minimizer in \eqref{eq:f1}.

\vspace{-10pt}
\subsection{Back from Gordon's Optimization to the LASSO}
\vp
Once we have analyzed (GO), we may now appropriately apply the Gaussian min-max Theorem (in particular, \cite[Thm. II.1]{tight}) to conclude with the following result.
\vspace{-3pt}
\begin{thm}[LASSO Objective]\label{thm:obj0}
 Recall the definitions of $\Lc(\w)$, $\Dc(\alpha,\beta,p)$ in \eqref{eq:main} and \eqref{eq:D}, respectively. Also, let
$\Lc_*:= \min_\w \Lc(\w)$ and $ \Dc_*:=\min_{\alpha\geq 0}\max_{\substack{0\leq\beta\leq 1\\ p>0}}\Dc(\alpha,\beta,p).$
Then, for all $\eps>0$,
$\lim_{n\rightarrow\infty}\Pro( |\Lc_*-\Dc_*|>\eps\Dc_* ) = 0.$
\end{thm}
\vspace{-6pt}
\noindent The proof follows from \eqref{eq:f1} when combining \eqref{eq:low} and \eqref{eq:up}. 
Also, a single application of \eqref{eq:low} in \eqref{eq:f2} shows that for any  $\R_\eps\subset\R$:
\begin{align}
\hspace{-5pt}\Lc_\eps: =\min_{\|\w\| \in \Rc_\eps} \Lc(\w) &\gtrsim \min_{\alpha\in\Rc_\eps}\max_{\substack{0\leq\beta\leq 1\\ p>0}}\Dc(\alpha,\beta,p) =:\Dc_\eps\label{eq:!2}.
\end{align}

To appreciate the power of Theorem \ref{thm:obj0}, note that it gives an explicit evaluation of the limiting behavior of the LASSO  optimal cost in terms of the optimal cost of a much simpler, scalar and deterministic optimization problem.
What is more, we can combine Theorem \ref{thm:obj0} with \eqref{eq:!2}, to prove the following result about the limiting behavior of the LASSO error.
\begin{thm}[LASSO error]\label{thm:err0}
Let $\hat\x$ be a minimizer of the LASSO in \eqref{eq:LASSO} and recall the definition of function $\Dc(\alpha,\beta,p)$ \eqref{eq:D}. If $\alpha_*$ is optimal for $\min_{\alpha\geq 0}\max_{\substack{0\leq\beta\leq 1\\ p>0}}\Dc(\alpha,\beta,p)$, then for all $\eps>0$:
$\lim_{n\rightarrow\infty}\Pro( |\|\hat\x-\x_0\|_2 - \alpha_*|>\eps\alpha_* ) = 0.
$
\end{thm}
\vspace{-5pt}
 We provide a sketch of the proof here (see also \cite[Thm.~II.1]{tight}).
Consider the event

$
\Ec= \left\{  \Lc_*\leq (1+\zeta_1)\Dc_* \text{ and } \Lc_\eps\geq (1-\zeta_2)\Dc_\eps \right\}
$.

\noindent In view of Thm.~\ref{thm:obj0} and \eqref{eq:!2}:
   $\lim_{n\rightarrow\infty}\Pro(\Ec) = 1$. Hence,
\vspace{-10pt}
\begin{align}
&\lim_{n\rightarrow\infty}\Pro\left( \Lc_\eps \leq (1+\delta ) \Lc_* \right) \leq\lim_{n\rightarrow\infty} \Pro\left( \Lc_\eps \leq (1+\delta ) \Lc_*  |\Ec\right)\nn\\
&\quad\quad\leq \Pro\left( ~(1-\zeta_2)\Dc_\eps \leq (1+\delta)(1+\zeta_1) \Dc_*~ \right).\label{eq:ab44}
\end{align}
Next, we show that we can choose $\zeta_1,\zeta_2,\delta$ such that the (deterministic) event in \eqref{eq:ab44} does not occur; this will complete the proof of \eqref{eq:prob_state}.
It can be shown 
 that the function $\Dc(\cdot,\beta,p)$ is strictly convex for all $\beta,p$. Thus, $\exists$ constant $L>0$ :
 \vspace{-10pt}
\begin{align}
&\Dc_\eps - \Dc_* =: \Dc(\alpha_\eps,\beta_\eps,p_\eps) - \Dc(\alphas,\betas,p_*) = \nn\\
&{\Dc(\alpha_\eps,\beta_\eps,p_\eps) - \Dc(\alpha_\eps,\beta_*,p_*)} + {\Dc(\alpha_\eps,\beta_*,p_*)- \Dc(\alphas,\betas,p_*)}\nn\\
&\vspace{-5pt}\geq 0 + L|\alpha_\eps - \alpha_*|^2 \geq \eps^2 L \alpha_*^2.\label{eq:str}
\end{align}
\vspace{-2pt}
Denote $\theta=\theta(\eps):=\eps^2 L \alpha_*^2/\Dc_*$. Set $\zeta_1=\theta/4$, $\zeta_2=\frac{\theta/2}{1+\theta}$, $\delta=\frac{\theta}{4+\theta}$ and apply \eqref{eq:str} to complete the proof.

\vspace{-10pt}
\subsection{Simplifying the Result}
\vspace{-7pt}
For any values of $\la$ and $\sigma$, we ask for an accurate prediction of the LASSO error $\|\hat\x^{\la,\sigma}-\x_0\|$. 
Theorem \ref{thm:err0} yields an answer
as the solution of a minimax optimization, which is scalar, deterministic and convex.
We have put some additional effort in order to simplify this optimization and make it more explicit. In particular, we can show that it has a unique global optimal solution, the evaluation of which essentially breaks down to solving two nonlinear one-dimensional algebraic equations (see Theorem \ref{thm:main}). Apart from the computational advantage of this alternative description, it also allows for further theoretical insights. For instance, starting from Theorem \ref{thm:main} we believe that it is possible to prove that the worst-case NSE is attained in the limit of the noise variance $\sigma^2\rightarrow 0$.


\vspace{-10pt}
\section{Results}

\vspace{-8pt}
\subsection{Arbitrary SNR Values}
\vspace{-7pt}
For the statement of our main result, recall \eqref{eq:D} and further consider the definitions in \eqref{eq:defn}. We also need the following.
\vspace{-9pt}
\begin{defn}[$\lac^\sigma$]\label{def:lac}
Define $\lac^\sigma$ as the unique solution of the equation $f(x,\frac{x}{\psi(x)})=0$. Further, define  $\lac^{0} := Q^{-1}\left(\frac{1}{2}\frac{m-k}{n-k}\right)$ and
\vspace{-14pt}
\begin{align*}
\lac^{min} := \begin{cases}
Q^{-1}\left(\frac{1}{2}\frac{m}{n-k}\right) & \text{ if } m<n-k,\\
0 & \text{else}.
\end{cases}
\end{align*}
 \end{defn}

\vspace{-10pt}
\noindent It can be shown that for all $\sigma>0$: $\lac^{min}\leq\lac^{\sigma}\leq\lac^{0}$.

\vspace{-5pt}
\begin{thm}\label{thm:main}
 Let $\hat\x^{\la,\sigma}$ be the solution of the LASSO in \eqref{eq:LASSO}.
  Recall the definitions of $f$ and $\lac^\sigma$ in \eqref{eq:defn} and Definition \ref{def:lac}.  Let $\lah := \max\{\la,\lac^\sigma\}$ and $q_*(\alpha) = \sqrt{1+\frac{m}{\sigma^2(m+\alpha^2)}}$. Denote $\alpha_*$ the unique solution of the equation $f\left(\lah, q_*(\alpha)\right)=0$ with respect to $\alpha$. Then,  the following limit holds in probability:
$\lim_{n\rightarrow\infty}\frac{\|\hat\x^{\la,\sigma}-\x_0\|_2}{\sigma\sqrt{m}}=\alpha_*.$
\end{thm}

\vspace{-20pt}
\subsection{Asymptotic NSE}\label{sec:as}
\vspace{-6pt}
Here we consider the case $\sigma^2\rightarrow 0$. 
\vp
\begin{thm}\label{thm:main2}
Let $\lah := \max\{\la,\lac^0\}$. 
If  $m\geq \Dc(\lah)$, then the following limit holds in probability
\vspace{-14pt}
$$\lim_{n\rightarrow\infty}\lim_{\sigma\rightarrow 0}\frac{\|\hat\x^{\la,0}-\x_0\|_2}{\sigma\sqrt{m}} = \sqrt\frac{{\Db(\lah)}}{{m-\Db(\lah)}} .$$
\end{thm}
\vspace{-7pt}
When $\la>\lac$, Theorem \ref{thm:main2} recovers the result of \cite[Theorem 3.2]{OTH};
adding to this, we are  able to characterize the behavior when $\la<\lac$.
 The theorem suggests that stable recovery is possible only when $m>\Dc(\lah)$. In particular, we require that at least $m>\min_{\la>0}\Dlf$. A recent line of work \cite{Sto,Cha,TroppEdge,stojnic2013upper}, has shown that $\min_{\la>0}\Dlf$ precisely characterizes the minimum number $m$ of required measurements for exact recovery of sparse signals under \emph{noiseless} linear measurements. Also, the minimum of the formula that appears in the theorem is achieved at $\la_{best}:=\arg\min_{\la>0}\Dlf$. We refer the reader to Figure \ref{fig:main} for an illustration of $\lac^0,\labb$ and to  \cite[Section 4.2]{OTH} for a detailed further discussion.


\vspace{-10pt}
\section{Simulation Results and Conclusion}

\vspace{-15pt}
\begin{figure}[!h]
\centering
\begin{subfigure}{.5\textwidth}
  \centering
  \includegraphics[width=0.9\linewidth]{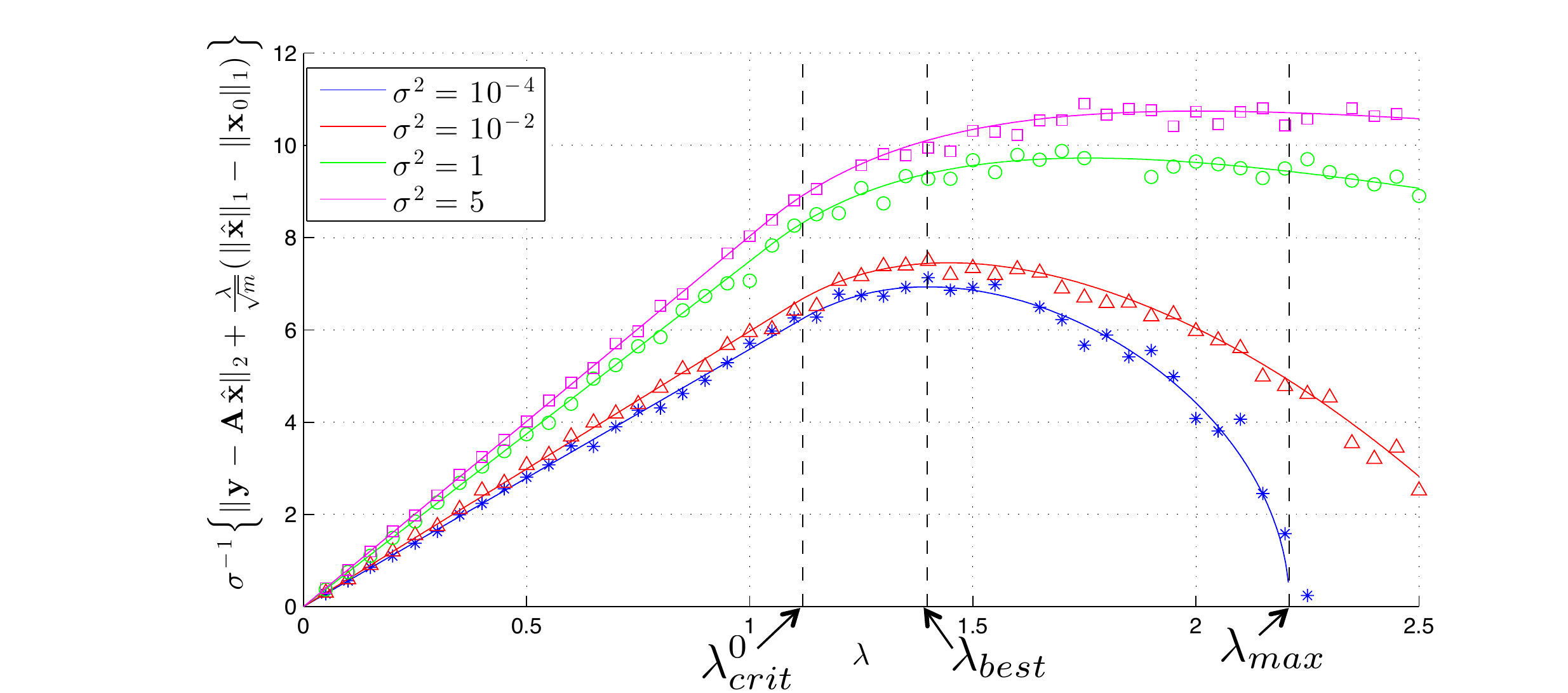}
  \caption{\footnotesize{optimal cost}}
  \label{fig:IIIA_1}
\end{subfigure}%
\\
\vspace{-2pt}
\begin{subfigure}{.5\textwidth}
  \centering
  \includegraphics[width=0.75\linewidth]{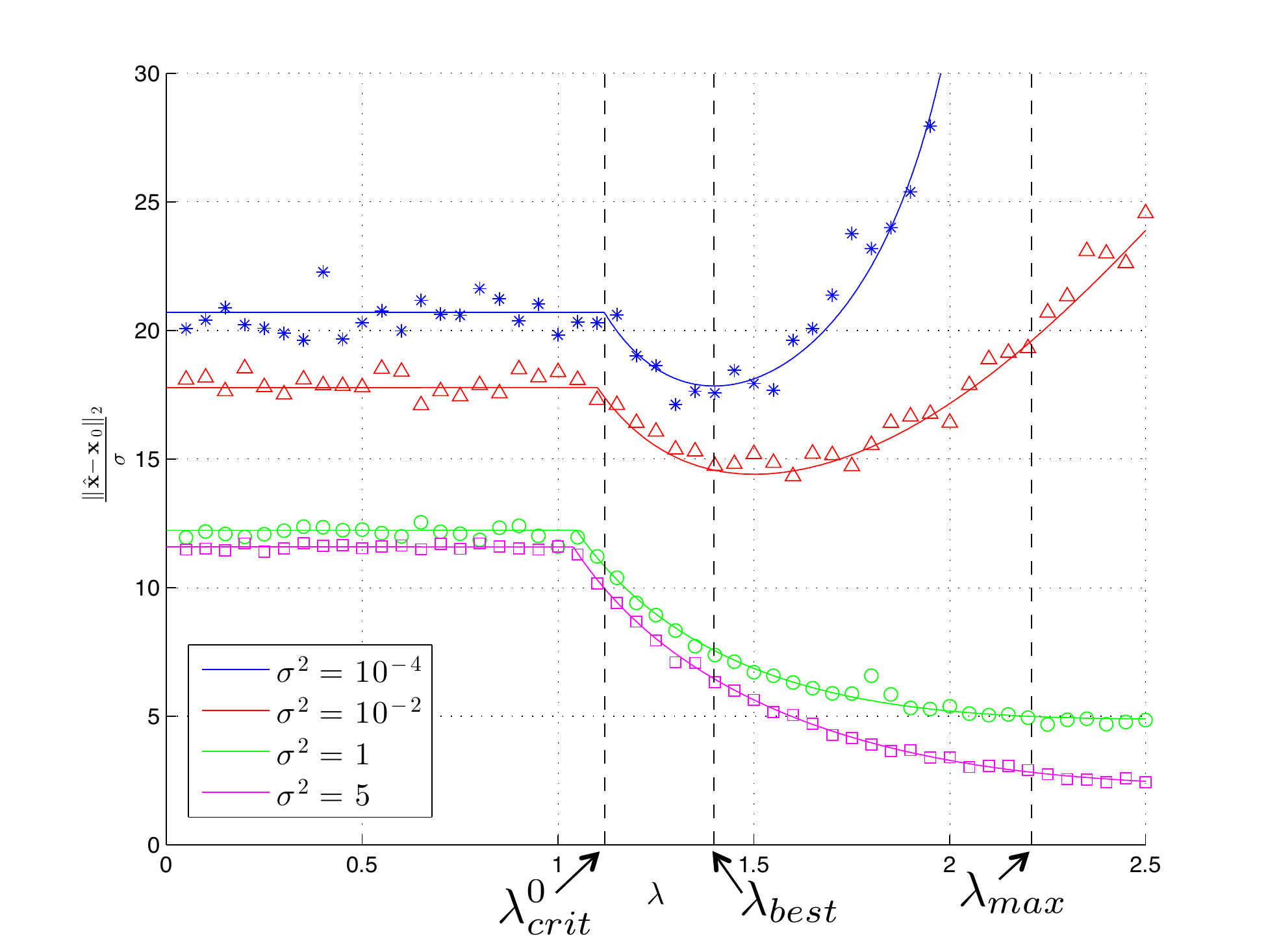}
  \caption{\footnotesize{normalized error}}
  \label{fig:IIIA_2}
\end{subfigure}
\vspace{-6pt}
  \caption{\footnotesize{$n=500$, $m = 150$, $k = 20$. Averages over 20 realizations.}}
\label{fig:main}
\end{figure}
\vspace{-10pt}
In  Figure \ref{fig:main}, observe the close agreement of the simulation results to the  predictions of Theorems \ref{thm:main} and \ref{thm:main2} (for the case $\sigma^2=10^{-4}$, we used Theorem \ref{thm:main2} for the prediction). Refer to Section \ref{sec:as} for the definitions of $\lac^0,\la_{best}$; $\la_{max}$ is such that $m<\Dlf$ for all $\la>\la_{max}$. Our expressions are accurate for problem sizes on the order of a few hundreds and valid for all values of $\sigma^2$. In Figure \ref{fig:IIIA_2}, the worst-case NSE occurs 
in the small noise-variance regime. We believe that this claim can be proved using Theorems \ref{thm:main} and \ref{thm:main2}.
\vp





\newpage
\bibliography{compbib}

\end{document}